\documentclass[12pt]{amsproc}
\usepackage{amssymb}
\usepackage{amscd}
\usepackage{amssymb,amsmath}
\usepackage{mathrsfs}

\newtheorem{definition}{Definition}[section]
\newtheorem{theorem}{Theorem}[section]
\newtheorem{lemma}{Lemma}[section]

\newtheorem{proposition}{Proposition}[section]

\def\bbi{\mathbb{I}}

\def\scr#1{\mathscr{#1}}

\def\ca{\mathscr{A}}

\def\cc{\mathscr{C}}

\def\ci{\mathscr{I}}

\def\cx{\mathscr{X}}

\def\Borel{\mbox{\rm Borel}}

\def\iff{\longleftrightarrow}
\def\then{\longrightarrow}

\begin{document}

\title{Completely nonmeasurable unions}

\author{Robert Ra\l owski and  Szymon \.Zeberski}

\address{ Institute of 
         Mathematics and Computer Science, 
         Wroc{\l}aw University of Technology, Wybrze\.ze Wyspia\'n\-skie\-go 27, 
         50-370 Wroc{\l}aw, Poland.}

\email[Robert Ra{\l}owski]{robert.ralowski@pwr.wroc.pl}
\email[Szymon \.Zeberski]{szymon.zeberski@pwr.wroc.pl}


\maketitle

\begin{abstract}
Assume that there is no quasi-measurable cardinal smaller than $2^\omega$.
($\kappa$ is quasi measurable if there exists $\kappa $-additive ideal
$\ci $ of subsets of $\kappa $ such that the Boolean
algebra $P(\kappa )/\ci$ satisfies c.c.c.)
 We show that for a c.c.c. $\sigma $-ideal $\bbi$ with a
Borel base of subsets of an uncountable Polish space, if
$\ca$ is a point-finite family of subsets from $\bbi$ then
there is an uncountable collection of pairwise disjoint subfamilies of $\ca$ whose union is completely
nonmeasurable i.e. its intersection with every non-small Borel set does not
belong to the $\sigma $-field generated by Borel sets and the
ideal $\bbi.$ This result is a generalization of Four Poles Theorem
(see \cite{BCGR}) and results from \cite{CMRRZ} and \cite{Z}.
\end{abstract}

\section{Notation and motivation}
In this paper
$X$ will denote an uncountable Polish space.
$\Borel$ will denote all Borel subsets of $X.$
A family $\bbi\subseteq P(X)$ will be a  $\sigma$-ideal of subsets of $X$ with Borel base containing singletons.
Let us recall that $\bbi $ has Borel base means that $(\forall I\in\bbi)(\exists J\in\bbi\cap\Borel)(I\subseteq J).$
We have the following cardinal coefficients
$$
\begin{array}{r@{\; =\; }l}
  \mbox{\rm add}(\bbi)   & \min\{|\cc|:\;\;\cc\subseteq\bbi,\;\bigcup\cc\notin\bbi\},    \\
  \mbox{\rm cov}(\bbi)   & \min\{|\cc|:\;\;\cc\subseteq\bbi,\;\bigcup\cc=X\},            \\
  \mbox{\rm cov}_h(\bbi) & \min\{|\cc|:\;\;\cc\subseteq\bbi,\;(\exists             
                 B\in\Borel\setminus\bbi)(\bigcup\cc\supseteq B)\},            \\
  \mbox{\rm cof}(\bbi)   & \min\{|\cc|:\;\;\cc\subseteq\bbi,(\forall I\in\bbi)(\exists 
                 C\in\cc)(I\subseteq C)\}.                                     \\
\end{array}
$$
Similarly  for a family $\ca\subseteq P(X)$ we can define
$$ 
\begin{array}{r@{\; =\; }l}     
   \mbox{\rm add}(\ca)        & \min\{|\cc|:\;\;\cc\subseteq\ca,\;\bigcup\cc\notin\bbi\}, \\
   \mbox{\rm cov}_h^\bbi(\ca) & \min\{|\cc|:\;\;\cc\subseteq\ca,\;(\exists        
                         B\in\Borel\setminus\bbi)(\bigcup\cc\supseteq B)\}.                \\
\end{array}
$$ 

We start our consideration with the following theorem from \cite{BCGR}.
It is known in literature as Four Poles Theorem.
\begin{theorem}[Brzuchowski, Cicho{\'n}, Grzegorek, Ryll-Nardzewski]
Let $\ca\subseteq\bbi$ be a point finite cover of $X$ i.e. $(\forall x\in X) |\{A\in\ca:\;\;x\in A\}|<\omega.$
Then there exists a subfamily $\ca'$ such that $\bigcup\ca'$ is not $\bbi$-measurable i.e. does not belong to 
the $\sigma$-field generating  by $\Borel$ and $\bbi.$
\end{theorem}

There is a hypothesis stated by J. Cicho{\'n} saying that we can improve the conclusion of the above theorem 
to get $\bigcup\ca'$ completely $\bbi$-nonmeasurable.

\begin{definition}
We say that $C$ is completely $\bbi$-nonmeasurable in $D$ iff
$$(\forall B\in\Borel\setminus\bbi)\left(B\cap D\notin\bbi\then (B\cap C\notin\bbi\wedge B\cap(D\setminus C)\notin\bbi)\right).$$ 
\end{definition}

Recall that $\bbi$ is c.c.c. if  every family $\ca\subseteq\Borel\setminus\bbi$ such that 
$$(\forall A,A'\in\ca)( A=A'\vee A\cap A'\in\bbi)$$
is at most countable.

If $\bbi$ is c.c.c then we can define
$[D]_\bbi$ to be a minimal (modulo $\bbi$) Borel set $B$ containing $D$ i.e.
$D\setminus B\in\bbi$ and if $D\subseteq C$ and $C$ is Borel then $B\setminus C\in\bbi.$

Assume that $\ca\subseteq\bbi.$ Let $\ci$ be an ideal on $P(\ca)$ associated with $\bbi$ in  the following way
$$
(\forall \cx\in P(\ca)) (\cx\in \ci\iff \bigcup \cx\in\bbi).
$$
Then $W\subseteq P(\ca)$ is an antichain in $P(\ca)/\ci$ iff $(\forall a,b\in W)(a\ne b\then a\cap b\in \ci)$. 
We say that $P(\ca)/\ci$ is c.c.c. iff every antichain on $P(\ca)/\ci$ is at most countable.

We say that the cardinal number $\kappa $ is {\it
quasi-measurable}  if there exists $\kappa $-additive ideal
$\ci$ of subsets of $\kappa $ such that the Boolean
algebra $P(\kappa )/\ci$ satisfies c.c.c. Cardinal $\kappa
$ is {\it weakly inaccessible} if $\kappa $ is regular cardinal
and for every cardinal $\lambda <\kappa$ we have that $\lambda^+
<\kappa.$ Recall that every quasi-measurable cardinal is weakly
inaccessible (see \cite{J}), so it is a large cardinal.

Let us recall a result from \cite{Z}.
\begin{theorem}[{\.Z}eberski]\label{zeberski}
Assume that there is no quasi-measurable cardinal not greater than
$2^\omega.$ Assume that $\bbi$ satisfies
c.c.c. Let $\ca\subseteq\bbi$ be a point-finite family such
that $\bigcup\ca\notin\bbi.$ Then there exists a subfamily
$\ca'\subseteq\ca$ such that $\bigcup\ca'$
is completely $\bbi$-nonmeasurable in
$\bigcup\ca.$
\end{theorem}

\section{Results}

Let us recall three technical lemmas from \cite{Z} (Theorem 3.3, Lemma 3.4, Lemma 3.5). 

\begin{lemma}[\.Zeberski]\label{bool}
Assume
that $\bbi$ satisfies c.c.c. Let $\{A_\xi
:\xi\in\omega_1\}$ be any family of subsets of $X.$ Then we can
find a family $\{I_\alpha\}_{\alpha\in\omega_1} $ of pairwise
disjoint countable subsets of $\omega_1$ such that for $\alpha
<\beta <\omega_1$ we have that $[\bigcup_{\xi\in I_\alpha
}A_\xi]_{\bbi}=[\bigcup_{\xi\in I_\beta }A_\xi]_{\bbi}.$
\end{lemma}

Next lemma is a reformulation of a result obtained in \cite{Z}.

\begin{lemma}[\.Zeberski]\label{otoczka}
Assume that $\bbi$ satisfies
c.c.c.  Let $\ca\subseteq\bbi$ be a point-finite family
such that $\bigcup\ca\notin\bbi$ and the algebra $P(\ca)/\ci $ is not c.c.c. Then there exists a
family $\{\ca_\alpha\}_{\alpha\in\omega_1}$ satisfying the
following conditions
\begin{enumerate}
\item $(\forall\alpha <\omega_1
)(\ca_\alpha\subseteq\ca\wedge
\bigcup\ca_\alpha\notin\bbi ),$
\item $(\forall\alpha <\beta <\omega_1)(\ca_\alpha\cap
\ca_\beta =\emptyset), $
\item $(\forall\alpha,\beta <\omega_1)([\bigcup\ca_\alpha ]_\bbi
=[\bigcup\ca_\beta ]_\bbi). $
\end{enumerate}
\end{lemma}

\begin{lemma}[\.Zeberski]\label{cienki}
Assume that $\bbi$ satisfies
c.c.c. Let $\ca\subseteq P(X)$ be any point-finite family.
 Then there exists a subfamily
$\ca'\subseteq\ca$ such that
$|\ca\setminus\ca'|\le\omega $ and $$(\forall B\in
\Borel\setminus\bbi )(\forall A\in\ca')
(B\cap\bigcup\ca\notin\bbi\rightarrow\neg
(B\cap\bigcup\ca\subseteq B\cap A)).$$
\end{lemma}

In paper \cite{CMRRZ} it is shown that if $\mbox{\rm cov}_h(\bbi)=\mbox{\rm cof}(\bbi)$ and $\ca\subseteq\bbi$ is a cover of $X$ such that 
$\bigcup\{A\in\ca:\;\;x\in A\}\in\bbi$ for every $x\in X,$ then there is a family $\ca'\subseteq\ca$ such that 
$\bigcup\ca'$ is completely $\bbi$ nonmeasurable.
This result can be generalized. Namely, we have the following theorem.
\begin{theorem}\label{ralowski}
 Let $\ca\subseteq\bbi$ be a family satisfying the following conditions:
\begin{enumerate}
 \item $(\forall B\in\Borel\setminus\bbi) |\{\ca(x) : \bigcup\ca(x)\cap B\neq\emptyset, 
                                         x\in X\}|=2^\omega,$
 \item $\mbox{\rm cov}_h^\bbi(\{ \bigcup\ca(x):\;\; x\in X\})=2^\omega,$ 
\end{enumerate}
where $\ca(x)=\{ A\in\ca:\; x\in A\}.$
Then there exists continuum many pairwise disjoint subfamilies $\{\ca_\alpha:\;\;\alpha\in 2^\omega\}$ of a family $\ca$ 
such that for every $\alpha\in 2^\omega$ a set $\bigcup\ca_\alpha$ is completely $\bbi$-nonmeasurable.
\end{theorem}
\begin{proof} Let us enumerate the set of all Borel $\bbi$ positive sets $\Borel\setminus\bbi=\{ B_\alpha:\; \alpha<2^\omega\}.$ 
By transfinite induction we will construct a sequence 
$$
( (A_{\xi,\eta},d_\xi )\in\ca\times B_\xi :\;\;\xi,\eta<2^\omega )
$$
with the following conditions:
\begin{enumerate}
 \item $(\forall \xi,\eta<2^\omega)( A_{\xi,\eta}\cap B_\xi\ne\emptyset)$,
 \item $\bigcup_{\xi,\eta<2^\omega} A_{\xi,\eta}\cap \{ d_\xi:\xi<2^\omega\}=\emptyset$,
 \item $(\forall\xi,\xi'<2^\omega)(\forall\eta,\eta' <2^\omega)(\eta\ne\eta'\longrightarrow A_{\xi,\eta}\ne A_{\xi',\eta'})$.
\end{enumerate}

Let us fix $\alpha< 2^\omega$ and assume that we have defined the sequence
$$
( (A_{\xi,\eta},d_\xi)\in\scr{A}\times B_\xi:\;\;\xi,\eta<\alpha )
$$
with the following conditions:
\begin{enumerate}
 \item $(\forall \xi,\eta<\alpha)( A_{\xi,\eta}\cap B_\xi\ne\emptyset)$,
 \item $\bigcup_{\xi,\eta<\alpha} A_{\xi,\eta}\cap \{ d_\xi:\xi<\alpha\}=\emptyset$,
 \item $(\forall\xi,\xi'<\alpha)(\forall\eta,\eta' <\alpha)(\eta\ne\eta'\longrightarrow A_{\xi,\eta}\ne A_{\xi',\eta'})$.
\end{enumerate}
For every $\xi<\alpha$ let us consider the set $\ca(d_\xi)=\{ A\in\ca:\;\; d_\xi\in A\}$.
By assumption (2) the family $\bigcup_{\xi<\alpha}\ca(d_\xi)$ does not cover any $\bbi$-possitive Borel set. 
So, assumptions (1) and (2) guaranties that  we can choose the set $\{ A_{\alpha,\eta}\in\ca:\;\; \eta<\alpha\}$ such that
\begin{enumerate}
 \item $(\forall\xi,\xi'\le\alpha)(\forall\eta,\eta'<\alpha)(\eta\ne\eta'\then A_{\xi,\eta}\ne A_{\xi',\eta'}),$ 
 \item $(\forall\eta<\alpha)(A_{\alpha,\eta}\cap B_\alpha\ne\emptyset),$ 
 \item $(\forall\xi,\eta<\alpha)( d_\xi\notin A_{\alpha,\eta}).$
\end{enumerate}
The same argument gives us the set $\{ A_{\xi,\alpha}\in\scr{A}:\;\; \xi\le \alpha\}$ with the following properties:
\begin{enumerate}
 \item $(\forall \xi,\xi'\le \alpha)(\forall\eta<\alpha)(A_{\xi,\eta}\ne A_{\xi',\alpha})$,
 \item $(\forall \xi\le \alpha)( A_{\xi,\alpha}\cap B_\xi\ne\emptyset\;\land\; A_{\xi,\alpha}\cap\{ d_{\xi'}:\xi'<\alpha\}=\emptyset).$
\end{enumerate}
Once again by assumption (2) we can find $d_\alpha\in B_\alpha$ such that 
$(\bigcup_{\xi,\eta\le \alpha} A_{\xi,\eta})\cap \{ d_\alpha\}=\emptyset$. 
It finishes the $\alpha$-step of our construction.

Now, let us put $\ca_\eta=\{ A_{\xi,\eta}\in\ca:\;\; \xi< 2^\omega\}$ for any $\eta<2^\omega.$
The family $\{ \ca_\eta:\;\eta< 2^\omega\}$  fulfills the assertion of our Theorem.
\end{proof}

\begin{proposition}
If $\mbox{ \rm cov}_h(\bbi)=2^\omega$ is a regular cardinal and $\ca\subseteq\bbi$ is a cover of $X$ such that each point is covered by less 
than continuum many members of $\ca$ then 
there exists continuum many pairwise disjoint subfamilies $\{\ca_\alpha:\;\;\alpha\in 2^\omega\}$ of a family $\ca$ 
such that for every $\alpha\in 2^\omega$ a set $\bigcup\ca_\alpha$ is completely $\bbi$-nonmeasurable.
\end{proposition}

\begin{proposition}
If $\mbox{ \rm cov}_h(\bbi)=2^\omega$  and $\ca\subseteq\bbi$ is  a point-finite family such that $\bigcup\ca\notin\bbi$ then 
there exists continuum many pairwise disjoint subfamilies $\{\ca_\alpha:\;\;\alpha\in 2^\omega\}$ of a family $\ca$ 
such that for every $\alpha\in 2^\omega$ a set $\bigcup\ca_\alpha$ is completely $\bbi$-nonmeasurable in $\bigcup\ca.$
\end{proposition}
\begin{proof}
 First, we use Theorem \ref{ralowski} to obtain continuum many pairwise disjoint subfamilies $\ca_\alpha$ for $\alpha<2^\omega$ such that $\bigcup\ca_\alpha$ is completely $\bbi$-nonmeasurable in $[\bigcup\ca]_\bbi.$ Then by point-finiteness of  family $\ca$ the  family $\{\bigcup\ca_\alpha\}_{\alpha<2^\omega}$ is also point-finite.
Using Lemma \ref{cienki} we can find a countable set $C\in [2^\omega ]^\omega$ such that  each member of the family
$\{\bigcup\ca_\alpha:\;\alpha\in 2^\omega\setminus C\}$ does not contain any $\bbi$-possitive Borel set with respect to $\bigcup\ca.$ So, the family $\{\ca_\alpha:\;\alpha\in 2^\omega\setminus C\}$ satisfies required conditions.
\end{proof}

Recall that the $\sigma$-ideal $\bbi$ has Steinhaus property if for any two $\bbi$-positive Borel sets $A,B\in\Borel\setminus\bbi$ the complex sum $A+B=\{a+b:\;a\in A, b\in B\}$ contains nonempty open set.
Let us remark that if the ideal $\bbi$ has Steinhaus property then $\mbox{\rm cov}_h(\bbi)=\mbox{\rm cov}(\bbi).$ 

\begin{theorem}\label{niequasi}
 Assume that $2^\omega$ is the smallest quasi-measurable cardinal. Let $\ca\subseteq\bbi$ be a point-finite family such that $\bigcup\ca\notin\bbi.$ Then $P(\ca)/\ci$ is not c.c.c.
\end{theorem}
\begin{proof}
 Assume that $\ca\subseteq\bbi$ satisfies the following conditions
 \begin{enumerate}
  \item $\bigcup\ca\notin\bbi,$
  \item $P(\ca)/\ci$ is c.c.c. 
 \end{enumerate}
 Since $2^\omega$ is the minimal quasi-measurable cardinal,  $|\ca|=2^\omega.$ 
Moreover $\mbox{\rm add}(\ca)=2^\omega$. By point-finiteness of the family $\ca$ we get that 
$\mbox{\rm add}(\{\bigcup\ca(x):\;\; x\in X\})=2^\omega,$ where $\ca(x)=\{A\in\ca:\;\; x\in A\}.$
So the family $\ca$ fulfils the assumptions of Theorem \ref{ralowski} (for $X=[\bigcup\ca]_\bbi$).
By Theorem \ref{ralowski} there exists $\{ \cc_\alpha:\; \alpha< 2^\omega\}$ such that
\begin{enumerate}
 \item $\cc_\alpha\subseteq \ca$ for any $\alpha<2^\omega$,
 \item $\forall\alpha<2^\omega\;\; \bigcup\cc_\alpha$ is completely $\bbi$-nonmeasurable in $[\bigcup\ca]_\bbi$,
 \item $\forall \alpha,\beta<2^\omega\;\;\alpha\ne\beta\then \cc_\alpha\cap\cc_\beta=\emptyset$.
\end{enumerate}
In particular, a family $ \{ \cc_\alpha:\; \alpha< 2^\omega\}$ forms an antichain in $P(\ca)/\ci,$ what gives a contradiction.
\end{proof}

\begin{theorem}\label{main}
 Assume there is no quasi-measurable cardinal smaller than $2^\omega$. Assume that the ideal $\bbi$ is c.c.c. Let  $\ca\subseteq\bbi$ be a family satisfying the following conditions:
\begin{enumerate}
 \item $\bigcup\ca\notin\bbi,$
 \item $(\forall x\in X) |\{A\in\ca:\;\;x\in A\}|<\omega.$
\end{enumerate}
Then there exists pairwise disjoint subfamilies $\{ \ca_\xi:\xi\in\omega_1\}$ of a family $\ca$ such that each of the union $\bigcup\ca_\xi$ is completely $\bbi$-nonmeasurable in $\bigcup\ca.$
\end{theorem}

\begin{proof} 
By transfinite induction we construct a family $\{B_\alpha\}$ of
pairwise disjoint Borel sets and a family
$\{\{\ca_\xi^\alpha\}_{\xi\in\omega_1 }\}$ of  subfamilies
of $\ca$ satisfying the following conditions
\begin{enumerate}
\item $B_\alpha\cap\bigcup\ca\notin\bbi,$
\item $(\forall\xi <\zeta
<\omega_1)(\ca_\xi^\alpha\cap\ca_\zeta^\alpha=\emptyset),$
\item $(\forall\xi <\omega_1)([\bigcup\ca^\alpha_\xi\setminus\bigcup_{\beta
<\alpha}B_\beta]_\bbi=B_\alpha).$
\end{enumerate}
At $\alpha $-step we consider the family
$\ca^\alpha=\{A\setminus\bigcup_{\xi <\alpha }B_\xi :
A\in\ca\setminus\bigcup_{\xi<\alpha}\ca_\xi\}.$ 
If $\bigcup\ca^\alpha\in \bbi$ then we
finish our construction. 
If $\bigcup\ca^\alpha\notin\bbi$
then by Theorem \ref{niequasi} the algebra $P(\ca^\alpha)/\ci$ is not c.c.c. We use Lemma \ref{otoczka} to obtain a required family
$\{\ca_\xi^\alpha\}_{\xi \in\omega_1 }.$ We put
$B_\alpha=[\bigcup\ca_0^\alpha\setminus\bigcup_{\zeta
<\alpha }B_\zeta]_\bbi.$

Since $\bbi$ satisfies c.c.c. the construction
have to end up at some step $\gamma <\omega_1.$

Now put $\ca'_\xi=\bigcup_{\alpha<\gamma}\ca^\alpha_\xi.$
By construction for each $\xi<\omega_1$ we have $[\bigcup \ca'_\xi]_\bbi=\bigcup_{\alpha<\gamma}B_\alpha=[\bigcup\ca]_\bbi.$

The family $\{\bigcup\ca'_\xi:\xi\in\omega_1\}$ is point finite because for every $x\in X$
$$
\left|\left\{\bigcup\ca'_\xi: x\in \bigcup\ca'_\xi\right\}\right|\le |\{A\in\ca:x\in A\}|<\omega.
$$
Now using Lemma \ref{cienki} we can find a countable set 
 $C\in [\omega_1 ]^\omega$ such that  each member of the family
$\{\bigcup\ca'_\xi:\;\xi\in \omega_1\setminus C\}$ does not contain any $\bbi$-positive set of the form $B\cap \bigcup\ca,$ where $B$ is Borel. So, the family $\{\ca'_\xi:\;\xi\in \omega_1\setminus C\}$ satisfies required conditions.
\end{proof}


\begin{thebibliography}{9}

\bibitem{BCGR}  {\sc J. Brzuchowski, J. Cicho\'{n}, E. Grzegorek, C. Ryll-Nardzewski},
{\it On the existence of nonmeasurable unions}, Bull. Polish Acad. Sci. Math.
{\bf 27} (1979), 447-448.

\bibitem{CMRRZ}  {\sc J. Cicho\'{n}, M. Morayne, R. Ra\l{}owski, C.
Ryll-Nardzewski, S. \.Zeberski},
       {\it On nonmeasurable unions, } Topol. and Its Appl.
       {\bf 154} (2007), 884-893.

\bibitem{J} {\sc T. Jech}, {\it Set Theory, } Springer-Verlag,
(2003).

\bibitem{Z} {\sc Sz. \.Zeberski}, {\it On completely nonmeasurable unions}, Math. Log. Quart. {\bf 53} No. 1 (2007), 38-42. 

\end{thebibliography}
\end{document}